# ON TYPE III GENERALIZED HALF LOGISTIC DISTRIBUTION

A.K. Olapade

Department of Mathematics, Obafemi Awolowo University,

Ile-Ife, NIGERIA.

e-mail: akolapad@oauife.edu.ng

**Abstract**

It is well known that generalized models is attracting the attention of researchers in recent times because of their flexibilities. Particularly, the logistic model has been generalized and applied by many authors while the half logistic distribution has not recieve much attention in term of its generalization. In this paper, we considered a generalized form of half logistic model called type III generalized half logistic distribution. We obtained its probability density function, the cumulative distribution function, the $n^{th}$ moment, the median, the mode, the $100p$-percentage points and the order statistics of the generalized distribution are established.

*2000 Mathematics Subject Classifications:* Primary 62E15; Secondary 62E10.

KEY WORDS: Half logistic distribution, Moments, Order statistics.

## 1 Introduction

One of the probability distributions which is a member of the family of logistic distribution is the half logistic distribution with probability density function.

$$f_Y(y) = \frac{2e^y}{(1+e^y)^2}, \quad 0 < y < \infty, \tag{1.1}$$

and cumulative distribution function

$$F_Y(y) = \frac{e^y - 1}{1 + e^y}, \quad 0 < y < \infty. \tag{1.2}$$

Balakrishnan (1985) studied order statistics from the half logistic distribution. Balakrishnan and Puthenpura (1986) obtained best linear unbiased estimator of location and scale parameters of the half logistic distribution through linear functions of order statistics. Balakrishnan and Wong (1991) obtained approximate maximum likelihood estimates for the location and scale parameters of the half logistic distribution with Type-II Right-Censoring. The hazard function for half logistic distribution is given as

$$h(y) = 1 - F(y) = 1 - (e^y - 1)/(e^y + 1) = 2/(e^y + 1). \tag{1.3}$$

Olapade (2003) proved some theorems that characterized the half logistic distribution. The half logistic distribution has not receive much attention from researchers in terms of generalization. Hence, in this research work, we present a new generalization of the half logistic distribution and obtain its properties, we also state and prove two theorems that characterize the generalized distribution.

## 2 Type III generalized half logistic distribution

Balakrishnan and Leung (1988) obtained three generalized forms of the standard logistic distribution one of which is presented below:

$$f_Y(y; b) = \frac{1}{B(b,b)} \frac{e^{by}}{(1 + e^y)^{2b}}, \quad -\infty < y < \infty, \; b > 0, \tag{2.1}$$

The probability density function above is called type III generalized logistic distribution. In this research, we want to obtain a generalized half logistic form of the generalized logistic distribution in equation (1.1) which we shall call type III generalized half logistic distribution and study its theories and properties.

The above equation is a symmetrical distribution along the real line and like the standard logistic distribution which is also symmetrical can be folded into two to obtain the standard half logistic distribution, the above type III generalized logistic distribution in equation (2.1) is folded into two to obtain

$$f_X(x; b) = \frac{2}{B(b,b)} \frac{e^{bx}}{(1 + e^x)^{2b}}, \quad 0 < x < \infty, \; b > 0. \tag{2.2}$$

We refer to the probability density in the equation (2.2) above as the type III generalized half logistic distribution. As expected, equation (2.2) reduces to the standard half logistic distribution in equation (1.1) when $b = 1$.

## 2.1 The cumulative distribution function of the type III generalized half logistic distribution

The cumulative distribution function of the type III generalized half logistic distribution is obtained as

$$F_X(x) = \int_0^x f_X(t)dt = \frac{2}{B(b,b)} \int_0^x \frac{e^{bt}}{(1+e^t)^{2b}}dt. \tag{2.3}$$

The integral above is not close, hence it is obtained numerically. As the values of $F_X(x)$ depends on the values of $b$ and $x$, some values of $F_X(x)$ are tabulated in the table below.

Table 2: The table of cumulative distribution of the type III generalized half logistic distribution.

$b = 2$

| x   | 0.0    | 0.1    | 0.2    | 0.3    | 0.4    | 0.5    | 0.6    | 0.7    | 0.8    | 0.9    |
|-----|--------|--------|--------|--------|--------|--------|--------|--------|--------|--------|
| 0   | 0.0000 | 0.0749 | 0.1490 | 0.2217 | 0.2922 | 0.3600 | 0.4246 | 0.4855 | 0.5425 | 0.5953 |
| 1.0 | 0.6439 | 0.6881 | 0.7281 | 0.7641 | 0.7962 | 0.8246 | 0.8496 | 0.8716 | 0.8907 | 0.9072 |
| 2.0 | 0.9215 | 0.9338 | 0.9443 | 0.9532 | 0.9608 | 0.9672 | 0.9726 | 0.9772 | 0.9810 | 0.9842 |
| 3.0 | 0.9869 | 0.9892 | 0.9910 | 0.9926 | 0.9939 | 0.9950 | 0.9958 | 0.9966 | 0.9972 | 0.9977 |
| 4.0 | 0.9981 | 0.9984 | 0.9987 | 0.9989 | 0.9991 | 0.9993 | 0.9994 | 0.9995 | 0.9996 | 0.9997 |
| 5.0 | 0.9997 | 0.9998 | 0.9998 | 0.9999 | 0.9999 | 0.9999 | 0.9999 | 0.9999 | 0.9999 | 1.0000 |

$b = 3$

| x   | 0.0    | 0.1    | 0.2    | 0.3    | 0.4    | 0.5    | 0.6    | 0.7    | 0.8    | 0.9    |
|-----|--------|--------|--------|--------|--------|--------|--------|--------|--------|--------|
| 0   | 0.0000 | 0.0935 | 0.1856 | 0.2751 | 0.3606 | 0.4412 | 0.5161 | 0.5847 | 0.6468 | 0.7022 |
| 1.0 | 0.7510 | 0.7935 | 0.8301 | 0.8612 | 0.8875 | 0.9094 | 0.9275 | 0.9423 | 0.9544 | 0.9641 |
| 2.0 | 0.9719 | 0.9781 | 0.9830 | 0.9869 | 0.9899 | 0.9922 | 0.9941 | 0.9955 | 0.9966 | 0.9974 |
| 3.0 | 0.9980 | 0.9985 | 0.9989 | 0.9992 | 0.9994 | 0.9995 | 0.9996 | 0.9997 | 0.9998 | 0.9998 |
| 4.0 | 0.9999 | 0.9999 | 0.9999 | 0.9999 | 0.9999 | 0.9999 | 1.0000 | 1.0000 | 1.0000 | 1.0000 |

The probability that a type III generalized half logistic random variable $X$ lies in an interval $(\alpha_1, \alpha_2)$ is given as

$$pr(\alpha_1 < X < \alpha_2) = F_X(\alpha_2) - F_X(\alpha_1)$$

$$\frac{2}{B(b,b)}[\int_0^{\alpha_2} \frac{e^{bt}}{(1+e^t)^{2b}}dt - \int_0^{\alpha_1} \frac{e^{bt}}{(1+e^t)^{2b}}dt] \qquad (2.3)$$

for any real value of $b$ and any given interval $\alpha_1 < \alpha_2$.

## 3 Moments of the Type III Generalized Half Logistic Distribution

Considering the Type III generalized half logistic distribution function $f_X(x; b)$ given in equation (2.2). The $n^{th}$ moment of $X$ is given as

$$E[X^n] = \frac{2}{B(b,b)} \int_0^\infty \frac{x^n e^{bx}}{(1+e^x)^{2b}} dx. \qquad (3.1)$$

The table three below shows a tabulated value of $E[X^n]$ for $b = 1.00, 2.00, ..., 10.00$. These values can be used to compute the mean, variance, skewness and kurtosis for the type III generalized half logistic distribution.

Table 3: The table of moment of the type III generalized half logistic distribution.

| b | $E[X]$ | $E[X^2]$ | $E[X^3]$ | $E[X^4]$ |
|---|---|---|---|---|
| 1.0000 | 1.3863 | 3.2899 | 10.8185 | 45.4576 |
| 2.0000 | 0.8863 | 1.2899 | 2.5008 | 5.9791 |
| 3.0000 | 0.6988 | 0.7899 | 1.1713 | 2.1095 |
| 4.0000 | 0.5946 | 0.5677 | 0.7054 | 1.0564 |
| 5.0000 | 0.5263 | 0.4427 | 0.4825 | 0.6307 |
| 6.0000 | 0.4771 | 0.3627 | 0.3562 | 0.4182 |
| 7.0000 | 0.4395 | 0.3071 | 0.2766 | 0.2973 |
| 8.0000 | 0.4095 | 0.2663 | 0.2228 | 0.2221 |
| 9.0000 | 0.3850 | 0.2351 | 0.1844 | 0.1722 |
| 10.0000 | 0.3644 | 0.2104 | 0.1559 | 0.1374 |

# 4 The median of the type III generalized half logistic distribution

Since the median of a probability density function is $f(x)$ is a point $x_m$ such that $\int_{-\infty}^{x_m} f(x)dx = \frac{1}{2}$; that is $F(x_m) = 1/2$.

For the type III generalized half logistic distribution

$$F_X(x) = \int_0^x f_X(t)dt = 1/2$$

$$\Rightarrow \frac{2}{B(b,b)} \int_0^x \frac{e^{bt}}{(1+e^t)^{2b}} dt = 1/2. \tag{4.1}$$

This integral can be evaluated numerically using iterative method in computer programmes. By using Maple V computer programme, we obtain the median of the type III generalized half logistic distribution for values of $b = 1, 2, ..., 5$.

Table 4: Table of median of the type III generalized half logistic distribution.

| b | Median |
|---|---|
| 1.00 | 1.09865 |
| 2.00 | 0.75475 |
| 3.00 | 0.57784 |
| 4.00 | 0.49445 |
| 5.00 | 0.43907 |

# 5 The $100p$-percentage point of the type III generalized half logistic distribution

For the type III generalized half logistic distribution the 100p-percentage point is obtained from the equation

$$F(x_p) = \frac{2}{B(b,b)} \int_0^{x_p} \frac{e^{bt}}{(1+e^t)^{2b}} dt = p. \tag{5.1}$$

Hence the 100p-percentage point of the type III generalized half logistic distribution can be obtained using iterative method as done for the median for specified $p$ and $b$.

# 6 The mode of the type III generalized half logistic distribution

For the type III generalized half logistic distribution the mode is obtained as follows

$$f_X(x;b) = \frac{2}{B(b,b)} \frac{e^{bx}}{(1+e^x)^{2b}}, \quad 0 < x < \infty, \ b > 0.$$

$$f'_X(x;b) = \frac{2}{B(b,b)} \left[ \frac{be^{bx}(1+e^x)^{2b} - 2be^x e^{bx}(1+e^x)^{2b-1}}{(1+e^x)^{4b}} \right]. \tag{6.1}$$

Equating the derivative to zero and solve for $x$ gives

$$x = \ln 1 = 0. \tag{6.1}$$

This is so because type III generalized half logistic distribution is obtained from a symmetric probability distribution.

# 7 Order statistics from type III generalized half logistic distribution

Let $X_1, X_2, ..., X_n$ be $n$ independently continuous random variables from the type III generalized half logistic distribution and let $X_{1:n} \leq X_{2:n} \leq .... \leq X_{n:n}$ be the corresponding order statistics. Let $F_{X_{r:n}}(x), \ r = 1, 2, ..., n$ be the cumulative distribution function of the $r^{th}$ order statistics $X_{r:n}$ and $f_{X_{r:n}}(x)$ denotes its probability density function. David (1970) gives the probability density function of $X_{r:n}$ as

$$f_{X_{r:n}}(x) = \frac{1}{B(r, n-r+1)} P^{r-1}(x)[1-P(x)]^{n-r} p(x). \tag{7.1}$$

For the type III generalized half logistic distribution with probability density function and cumulative distribution function given as

$$f(x) = p(x) = \frac{2}{B(b,b)} \frac{e^{bx}}{(1+e^x)^{2b}}, \quad 0 < x < \infty, b > 0$$

and

$$F_X(x) = P(x) = \frac{2}{B(b,b)} \int_0^x \frac{e^{bt}}{(1+e^t)^{2b}} dt$$

respectively. By substitution into equation (7.1), we have the probability density function of the $r^{th}$ order statistics from the type III generalized half logistic distribution as

$$f_{X_{r:n}}(x) = \frac{1}{B(r, n-r+1)} \left[\frac{2}{B(b,b)} \int_0^x \frac{e^{bt}}{(1+e^t)^{2b}} dt\right]^{r-1} \left[1 - \frac{2}{B(b,b)} \int_0^x \frac{e^{bt}}{(1+e^t)^{2b}} dt\right]^{n-r}$$

$$\times \frac{2}{B(b,b)} \frac{e^{bx}}{(1+e^x)^{2b}}. \tag{7.2}$$

Consider the probability density of the $r^{th}$ order statistics from the type III generalized half logistic distribution in equation (7.2). When $r = n$, then the probability density function of the maximum order statistic is

$$f_{X_{n:n}}(x) = n[\frac{2}{B(b,b)} \int_0^x \frac{e^{bt}}{(1+e^t)^{2b}} dt]^{n-1} \frac{2}{B(b,b)} \frac{e^{bx}}{(1+e^x)^{2b}}$$

$$= \frac{2n}{B(b,b)} [\frac{2}{B(b,b)} \int_0^x \frac{e^{bt}}{(1+e^t)^{2b}} dt]^{n-1} \frac{e^{bx}}{(1+e^x)^{2b}}. \tag{7.3}$$

Also, consider the probability density of the $r^{th}$ order statistics from the type III generalized half logistic distribution in equation (7.1). When $r = 1$, then the probability density function of the minimum order statistic is

$$f_{X_{1:n}}(x) = \frac{2n}{B(b,b)} [1 - \frac{2}{B(b,b)} \int_0^x \frac{e^{bt}}{(1+e^t)^{2b}} dt]^{n-1} \frac{e^{bx}}{(1+e^x)^{2b}}. \tag{7.4}$$